\documentclass[11pt]{amsart}
\usepackage{graphicx}
\usepackage{latexsym}
\usepackage{amsfonts}
\usepackage{amsthm}
\usepackage{amsmath}

\newtheorem{theorem}{Theorem}[section]
\newtheorem{lemma}[theorem]{Lemma}

\theoremstyle{definition}
\newtheorem{definition}[theorem]{Definition}

\newtheorem{example}[theorem]{Example}

\def\op#1{\operatorname{#1}}
\def\ring#1{\mathbb{#1}}
\def\S#1{\ring{S}_{\op{#1}}}

\def\rq{\text{'}}
\def\lq{\text{`}}

\begin{document}

%% Title, optional running title in brackets.
\title{What is motivic measure?}
\author{Thomas C. Hales}

\begin{abstract}
This article gives an exposition of the theory of arithmetic
motivic measure, as developed by J. Denef and F. Loeser.
\end{abstract}

\maketitle

\section{Preliminary Concepts}

There is much that is odd about motivic measure if it is judged by
measure theory in the sense of twentieth century analysis.  It
does not fit neatly with the tradition of measure in the style of
Hausdorff, Haar, and Lebesgue.  It is best to view motivic measure
as something new and different, and to recognize that when it
comes to motivic measure, the term `measure' is used loosely.

Motivic integration has been developing at a break-neck pace, ever
since Kontsevich gave the first lecture on the topic in 1995.
This article gives an exposition of the theory of arithmetic
motivic measure, as developed by J. Denef and F. Loeser.

Motivic measure will be easier to understand, once two of its
peculiarities are explained.  The first peculiarity is that the
measure is not real-valued. Rather, it takes values in a scissor
group.  An introductory section on scissor groups for polygons
will recall the basic facts about these groups. The second
peculiarity is that rather than a boolean algebra of measurable
sets, we work directly with the underlying boolean formulas that
define the sets.  The reasons for working directly with boolean
formulas will be described in a second introductory section.

After these two introductory remarks, we will describe `motivic
counting' in Section~\ref{sec:counting}.  Motivic counting is to
ordinary counting what motivic measure is to ordinary measure.
Motivic counting will lead into motivic measure.

\subsection{Scissor Groups for polygons}

Motivic volume is defined by a process that is similar to the
scissor-group construction of the area of polygons in the plane.
To draw out the similarities, let us recall the construction. It
determines the area of polygons without taking limits.

\begin{figure}[htb]
  \centering
  \includegraphics{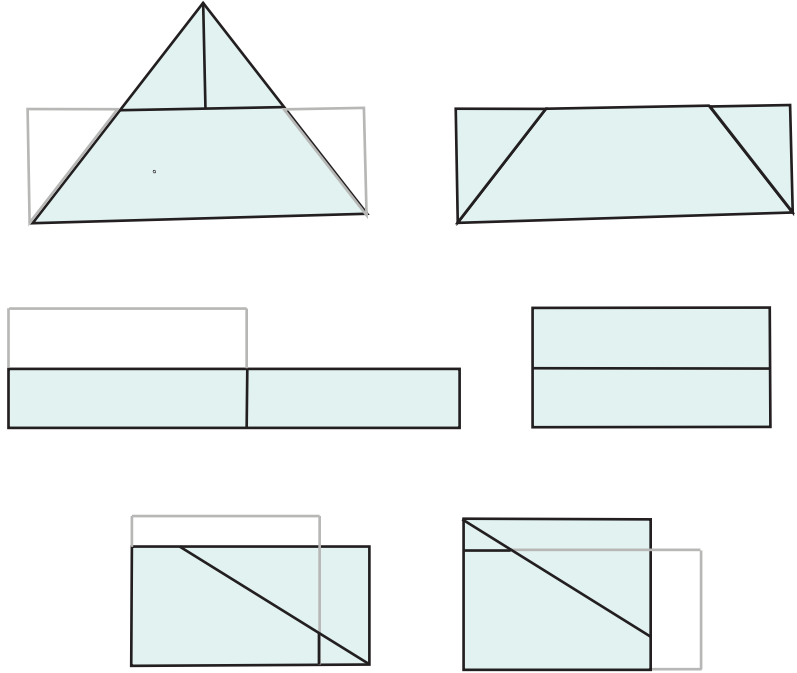}
  \caption{Triangles transform into
  unit width rectangles by scissor and congruence relations.  Later, we will transform
  ring formulas into algebraic varieties by scissor and congruence relations.}
  \label{fig:scissor}
\end{figure}

Any polygon in the plane can be cut into finitely many triangles
that can be reassembled into a rectangle of unit width.
Figure~\ref{fig:scissor} illustrates three steps ($2$, $3$, and
$4$) of the general algorithm. The algorithm consists of $5$
elementary transformations.  (1) Triangulate the polygon. (2)
Transform triangles into rectangles. (3) Fold long rectangles in
half. (4) Rescale each rectangle to give it an edge of unit width.
(5) Stack all the unit width rectangles end to end.  The length of
the unit width rectangle is the area.

An abelian group encodes these cut and paste operations. Let $F$
be the free abelian group on the set of polygons in the plane.
%For example, an element of $F$ is the formal sum
%    $$
%    3[\Box(1)] - 2[\triangle(3,4,5)],
%    $$
%three squares of edge length $1$ minus $2$ right triangles with
%legs $3$ and $4$.

%% Special Format:

We impose two families of relations:

\smallskip

\noindent {\bf Scissor relations}. If $P$ is a polygon that can be
cut into polygons $P_1$ and $P_2$, then
    $$[P] = [P_1]+[P_2]$$

\smallskip
\noindent {\bf Congruence relations}.  If $P$ and $P'$ are
congruent polygons then
    $$
    [P] = [P'].
    $$

\smallskip

The scissor group $\S{poly}$ of polygons is defined as the free
abelian group subject to these two families of relations.  In some
sense, this entire article is an exploration of scissor and
congruence relations in diverse contexts.  By and by, we will
construct several closely related scissor groups $\S{poly}$,
$\S{count}$, $\S{ring}$, $\S{cover}$, and $\S{mot}$, each
constructed as a free abelian group modulo scissor and congruence
relations.

\begin{theorem} The polygon scissor group $\S{poly}$
of polygons is isomorphic to the additive
group of real numbers $\ring{R}$.  Under this isomorphism, the
real number attached to the class $[P]$ of a polygon is its area.
\end{theorem}

\begin{proof}  A group homomorphism from
$\S{poly}$ to $\ring{R}$ sends each class $[P]$ to its area.  It
is onto, because there are polygons of every positive real area,
and negations of polygons of every negative real area.  By scissor
and congruence relations, every element of the scissor group is
represented by the difference of two unit width rectangles.  To be
in the kernel, the two rectangles must have the same area; but
then they are congruent, and their difference is the zero element
of $\S{poly}$. Thus, the homomorphism is also one-to-one.
\end{proof}

The area function on the set $\{P\}$ of polygons thus factors
through $\S{poly}$.
    \begin{equation}
        \begin{array}{lllll}
            \{P\}&\to &\S{poly}&\to&\ring{R}\\
            P &\mapsto &[P] &\mapsto &\op{area}(P)
        \end{array}
    \end{equation}

We might ponder which of these two maps ($P\mapsto[P]$ or
$[P]\mapsto\op{area}(P)$) captures the greater part of the
area-taking process.   Motivic measure commits to a position on
this issue: the first stage ($P\mapsto[P]$) is identified as the
area-taking process and the second stage $[P]\mapsto\op{area}(P)$
is a {\it specialization\/} of the area.  In this case,
specialization is an isomorphism.  Our approach to measure in this
article is decidedly unsophisticated: taking the measure of
something consists in mapping that thing into its scissor group,
$P\mapsto[P]$.

\subsection{The measure of a formula}

Traditionally, we take the measure of a set $X =\{X~|~\phi(x)\}$
(say a subset of a locally compact space), but we do not take the
measure of the formula $\phi$ defining a set.  With motivic
measure, we take the measure of the formula directly. Concretely,
the formula
    \begin{equation}
    \lq x^2+y^2=1 \rq
    \label{eqn:circle}
    \end{equation}
defines the circle
    \begin{equation}
    \{(x,y)~|~x^2+y^2=1\}.
    \label{eqn:circle-set}
    \end{equation}
With motivic measure, we take the measure of the equation of the
circle (Equation~\ref{eqn:circle}) rather than the measure of the
circle itself (Equation~\ref{eqn:circle-set}).  Attention shifts
from sets to formulas.

What purpose does it serve to measure formulas rather than the
underlying set?  As algebraic geometers are eager to remark, each
formula defines an infinite collection of sets. For instance, for
each finite field $\ring{F}_q$, we can take the set of
$\ring{F}_q$ points on the circle:
    \begin{equation}
    \{(x,y)\in\ring{F}_q^2~|~x^2+y^2=1\}.
    \end{equation}
We will see that the motivic measure of the formula is a universal
measure in the sense that the value it attaches to the formula
does not commit us to any particular field.   And yet if we are
supplied with a particular field, it will be possible to recover
the traditional measure of a set from the motivic measure of its
defining formula.  In this sense, motivic measure is to
traditional measures what an algebraic variety is to its set of
solutions.

\section{Counting measures and Finite Fields}\label{sec:counting}

Counting is the fountainhead of all measure.  The measure of a
finite set is its cardinality. At the risk of belaboring the
point, in preparation for what is to come, let us recast ordinary
counting.   The scissor relation for disjoint finite sets is
    $$
    [X\cup Y] = [X] + [Y].
    $$
More generally, if we allow the sets to intersect, it is
    \begin{equation}
    [X\cup Y] = [X] + [Y] - [X\cap Y].
    \label{eqn:scissor-set}
    \end{equation}
The congruence relation asserts that
    $$
    [X] = [X'].
    $$
whenever there is a bijection between $X$ and $X'$.  The scissor
group $\S{count}$ is the quotient of the free abelian group on
finite sets satisfying the scissor and congruence relations. It is
is isomorphic to $\ring{Z}$.  The cardinality $\#X$ of a finite
set $X$ factors through the scissor group
    $$
    X \mapsto [X] \mapsto \#[X]\in\ring{Z}.
    $$
Of course, if our only purpose were to count elements in finite
sets, this construction is overkill.   The first motivic measure
that we present is an analogue of this approach to counting. We
call it the {\it motivic counting measure}. The scissor relation
will be similar to Equation~\ref{eqn:scissor-set}.

\subsection{Ring formulas}

Traditional measure calls for a full discussion of the class of
measurable sets.  Since we work with formulas rather than sets,
our approach calls for a full discussion of the class of formulas
to be measured.

We allow all syntactically correct formulas built from a countable
collection of variables $x_i$, parentheses, and the symbols
    \begin{equation}
    \forall,\quad\exists,\quad\lor,\quad
    \land,\quad\lnot,
    \quad 0,\quad 1,\quad (+),\quad (-),\quad (*),\quad (=)
    \end{equation}
More precisely, we allow all formulas in the first-order language
of rings. %%
%% We are applying the special case of Denef and Loeser, where we
%% adjoin the constants of the field of rational numbers, but no more.
%%
A formula that has been constructed from these symbols will be
called a {\it ring formula}.  We avail ourselves of the usual
mathematical abbreviations and renamings of variables. We write
$3$ for $1+(1+1)$, $x^n$ for $x*x*x\cdots *x$ ($n$ times), $x y$
for $x*y$, $a+b+c$ for $a+(b+c)$, and so forth.

With usual abbreviations,
    $$
    \lq\forall x~y~z.~~x^3+y^3 = z^3\rq
    $$
is a ring formula, because its syntax is correct.  But
    $$
    \lq)) \forall +\forall = 2\forall ((\rq
    $$
and
    $$
    \lq\land\lor\land\lor\land\rq
    $$
are not ring formulas.

\subsection{The scissor group of ring formulas}

We imitate the construction of the scissor groups $\S{poly}$ and
$\S{count}$ to build the scissor group of ring formulas.

Take the free abelian group on the set of ring formulas.

We impose two families of relations.  The scissor relation takes
the form established in Equation~\ref{eqn:scissor-set} for unions.

\smallskip
\noindent{\bf Scissor relations}. If $\phi_1\lor\phi_2$ is a
disjunction of two formulas,  then
    \begin{equation}
    [\phi_1\lor\phi_2] = [\phi_1]+[\phi_2]-[\phi_1\land\phi_2].
    \end{equation}
    \smallskip

To describe the congruence relation, we must decide what it should
mean for two ring formulas to be congruent.  By way of analogy, in
the case of polygons, two are congruent if there is a bijection
between the two sets that is induced by an isometry.  Our first
guess at the congruence relation for ring formulas is that two
ring formulas are congruent if there is a bijection between the
sets of solutions for each finite field $\ring{F}_q$. (We limit
ourselves to {\it finite\/} fields because we are attempting to
imitate the counting measure of {\it finite\/} sets.)  However,
there are two modifications that we must make to this first guess
to arrive at a workable relation.

The first modification is to use pseudo-finite fields rather than
finite fields.  A {\it pseudo-finite field} is an infinite perfect
field such that every absolutely irreducible variety over the
field has a rational point and such that there is a unique field
extension of each finite degree (inside a fixed algebraic closure
of the field). The defining properties of a pseudo-finite field
are properties possessed by finite fields (except the part about
being infinite). Moreover, logicians have found that the behavior
of pseudo-finite fields is essentially no different from the
generic behavior of finite fields, but they avoid the hassles that
appear in positive characteristic.  For those seeing pseudo-finite
fields for the first time, it would not be a severe distortion of
the facts to ignore the `pseudo' and to work instead with finite
fields.

The second modification is to require the bijection between the
solutions to come from a ring formula that is independent of the
underlying field.  We are now ready to state the congruence
relations.

\smallskip \noindent{\bf Congruence relations}.
    $$[\phi] = [\phi']$$
if there exists a ring formula $\psi$ such that for every
pseudo-finite field $K$ of characteristic zero, the interpretation
of $\psi$ gives a bijection between the tuples in $K$ satisfying
$\phi$ and the tuples in $K$ satisfying $\phi'$.

\begin{example} The congruence relation gives
    $$
    [\lq\exists x.\quad x^2 + b x + c = 0\rq] = [\lq\exists X.\quad X^2
    = B^2 - 4 C\rq]
    $$
The formula $\psi$ realizing the congruence and the bijection at
the level of points is
    $$
    \lq(b=B)\land (c=C)\rq.
    $$
That is, in every pseudo-finite field of characteristic zero, a
monic quadratic polynomial has a root if and only if its
discriminant is a square.
\end{example}

\begin{definition}
The scissor group $\S{ring}$ of ring formulas is defined as the
free abelian group subject to the scissor and congruence
relations.
\end{definition}

\subsubsection{Counting measure}

\begin{definition}
The {\it counting measure\/} of a ring formula $\phi$ is its class
$[\phi]$ in the scissor group of ring formulas.
\end{definition}

\subsubsection{Fubini and Products}

There is a trivial sort of Fubini theorem for finite sets: the
cardinality of a Cartesian product of two sets is the product of
the cardinalities of the two sets.  To make sense of a Fubini
theorem for ring formulas, it is necessary to introduce products
to the scissor groups; that is, we need a {\it scissor ring}. This
is easy to arrange.  If $\phi_1(x_1,\ldots,x_n)$ is a formula with
free variables $x_1,\ldots,x_n$ and $\phi_2(y_1,\ldots,y_m)$ is a
formula with free variables $y_1,\ldots,y_m$, and if the free
variables of $\phi_1$ are distinct from the free variables of
$\phi_2$, then we declare the product to be
    $$
    \phi_1(x_1,\ldots,x_n)\land \phi_2(y_1,\ldots,y_m).
    $$
This induces a well-defined product\footnote{We have a moving
lemma: the congruence relation in the scissor group can be used to
relabel the free variables of a formula, so that free variables of
the two factors are always distinct.} on the scissor group
    \begin{equation}
    [\phi_1(x)][\phi_2(y)] = [\phi_1(x)\land\phi_2(y)].
    \label{eqn:fubini}
    \end{equation}
Under this product, the scissor group becomes a ring.
Equation~\ref{eqn:fubini} asserts that counting measure satisfies
a rather trivial Fubini theorem for ring formulas -- at least for
ring formulas without any shared free variables.

\subsubsection{The universal nature of the
counting measure}

The counting measure $[\phi]$ of a ring formula $\phi$ is designed
to be the {\it universal counting measure} for ring formulas.  For
every finite field $\ring{F}_q$, there is a special counting
measure on ring formulas:
    \begin{equation}
    \phi\mapsto \#_q(\phi) =
    \#\{(x_1,\ldots,x_n)\in \ring{F}_q^n~|~
        \phi^{\ring{F}_q}(x_1,\ldots,x_n)\}.
    \end{equation}
It gives the number of solutions to the ring formula over a
particular finite field.  In contrast, the general counting
measure of a ring formula takes values in a scissor ring whose
construction bundles all pseudo-finite fields together.

We can be precise about the way in which the counting measure
bundles the counting measures $\#_q(\phi)$. Each formula $\phi$
gives a function $q\mapsto \#_q(\phi)$, an integer-valued function
on the set of prime powers.  Let $F$ be the ring of all
integer-valued functions on the set $\{p^r\}$ of prime powers.
Declare two functions equivalent, if they take the same value at
$p^r$ for all $r$ and for all but finitely many $p$. Write
$F/\sim$ for the quotient of $F$ under this equivalence relation.

\begin{theorem}
There exists a ring homomorphism $N$ from the scissor ring
$\S{ring}$ to $F/\sim$ that respects counting: $\#_*(\phi) =
N([\phi])$.
\end{theorem}

In other words, with only a finite amount of ambiguity, the
counting measure specializes to counting solutions to ring
formulas over finite fields.  To say that $N$ is a ring
homomorphism is to say that it is compatible with products and
Fubini.  Unlike the earlier isomorphisms for polygons
$\S{poly}\cong \ring{R}$ and finite sets $\S{count}\cong\ring{Z}$,
here we make no claim of isomorphism between the scissor group
$\S{ring}$ and the target ring $F/\sim$.

The proof of the theorem relies on ultraproducts, a standard tool
in logic.

\subsection{Improving the scissor ring}
\label{sec:qe}

The shortcoming of the scissor ring $\S{ring}$ is that is too much
about it has been left inexplicit.  In our discussion of the area of
planar polygons, we found a handy set of generators (unit width
rectangles).  Our current aim is to find a handy set of generators
of a somewhat modified scissor ring $\S{mot}$.  The idea is to take
a ring formula, and through a process of ``quantifier elimination''
arrive at an equivalent ring formula that does not involve any
quantifiers (that is, the symbols $\forall$, $\exists$ will be
eliminated).  Quantifier elimination is a subject that was under
development long before motivic integration arrived on the scene.
Background on M. Fried and others' work on quantifier elimination
can be found in \cite{FAn} and \cite{FAr}, as well as in an appendix
to this article by M. Fried.

A formula without quantifiers belongs less to the realm of logic
than to the realm of algebraic geometry. A suggestive example of a
quantifier-free formula is
    $$
    (f_1=0)\land (f_2=0)\land \cdots \land (f_n=0).
    $$
That is, the zero set of an affine variety.  In fact, we will find
that the improved scissor ring is defined as a quotient of the
free abelian group on the set of varieties over $\ring{Q}$.  The
details of this construction will reveal what is so {\it
motivic\/} about motivic measure.

\subsection{A scissor ring for coverings}

%% This is slightly different from the covering relation used by Denef and
%% Loeser.  The relation used here is stronger (but easier to state).  To
%% see that you still get a homomorphism to $\S{mot}$, it is necessary
%% to revisit Theorem 3.3.5 of "Definable Sets" to check that if there is
%% an $n$-to-$1$ cover, then $\chi_c(A) = n \chi_c(B)$.  This is OK.

%% In more detail, 3.3.5 makes it so the n-1 correspondence psi is itself
%% a Galois stratification W/Z.  We can look stratum by stratum.  Going to larger covers
%% we may assume that the W'/Z' for phi and W/Z for psi is such that W/Z' is Galois.
%% There is a morphism W->Z->Z'.
%% Pull W' up to W so that the Galois covers are W/Z' and W/Z.
%% Use the ICM relation relating nonabelian and abelian covers so that W/Z is abelian with group C
%% and the W/Z relation is of the form phi[W,W/C,C] as in ICM.  In particular W/C = Z. (This scales the
%% constant n by an integer factor.)
%% The image of Z in Z' consists in elements that can be lifted at least to W/C.
%% If some point lifts to W/A with A sub C, then it is not in the image of phi[W,W/C,C], contrary to the
%% n-1 hypothesis.  Thus, the Z' formula is phi[W,Z',C].  The non-abelian/abelian ICM relation
%% rewrites phi[W,Z',C] in terms of phi[W,W/C,C], again changing n by an integer factor.
%% This shows that the cover relation is a consequence of the ICM-Denef-Loeser relations as desired.

Each ring homomorphism $f:\S{ring}\to R$ defines a specialization
of the counting measure
    $$
    \phi\mapsto[\phi]\to f[\phi] \in R.
    $$
The ring $F/\sim$ is one of many possible specializations $R$.

Another specialization of $\S{ring}$ comes from $n$-sheeted
covers:

\begin{definition} We say that one formula $\phi(x)$ is an
{\it $n$-sheeted cover\/} of another formula $\phi'(x')$ if there
exists a ring formula $\psi(x,x')$ such that for every
pseudo-finite field of characteristic zero, $\psi$ gives an $n$ to
$1$ correspondence between the solutions $x$ of $\phi(x)$ and the
solutions $x'$ of $\phi(x')$.
\end{definition}

\begin{example}  Let $\phi(x)$ be the formula $\lq x\ne 0\rq$ and let
$\phi'(y)$ be the formula
 $$
 \lq\exists~z.~ (z^2 = y) \land (y\ne 0)\rq.
 $$
The formula $\psi(x,y)$ given by
  $$
  \lq x^2 = y\rq,
  $$
presents $\phi$ as a $2$-sheeted cover of $\phi'$.
\end{example}

The congruence condition for $\S{ring}$ asserts that if $\phi$ is
a $1$-sheeted cover of $\phi'$, then they give the same class in
$\S{ring}$.  A broader congruence condition can be given as
follows.

\noindent{\bf Congruence (Covers).}  If $\phi$ is an $n$-sheeted
cover of $\phi'$ for some $n$, then
    $$
    [\phi] = n[\phi'].
    $$

We may form a new scissor ring $\S{cover}$ with this broader
congruence condition and the old scissor relation. We have a
canonical surjection $\S{ring}\to\S{cover}$.

\subsection{The scissor group of motives}

\hbox{}
\smallskip

\noindent{\bf Generators}.  Let $\op{Var}_{\ring{Q}}$ be the
category of varieties over the field of rational numbers
$\ring{Q}$.   We take the free abelian group generated by the
objects of $\op{Var}_{\ring{Q}}$.

An example of a element of the free abelian group is
$[\ring{A}^1]$, the generator attached to the affine line.  This
particular generator will be of special importance in the
constructions that follow.  We write $\ring{L}=[\ring{A}^1]$ for
this element and for its image in various scissor groups.  (The
`L' is for {\it Lefschetz}, as in Lefschetz motive.)

There are two types of relations: scissor relations and congruence
relations.  Our scissor relation will be rather crude, but
justifiably so, since the Zariski topology is a coarse topology
that limits the possibilities for a scissor relation.  The only
cutting that will be permitted is that of partitioning a variety
into a closed subvariety and its complement.

\smallskip\noindent
{\bf Scissor Relation}.  If $Z$ is a closed subvariety of $X$,
then
    $$
    [X] = [Z] + [X\setminus Z].
    $$
\smallskip

The congruence relation is more involved than the scissor
relation. If we make a direct translation of the congruence
relation for the scissor group of ring formulas, we might guess
that the congruence condition between two varieties $X$ and $Y$
should be the existence of a correspondence $\Psi$ between $X$ and
$Y$ that induces a bijection between $X(K)$ and $Y(K)$ for every
pseudo-finite field of characteristic zero. This first guess is
suggestive: the congruence relation should involve an algebraic
correspondence.  This suggestion lands us deep in the territory of
motives.  Here is the precise definition of the congruence
relation.

\smallskip
\noindent{\bf Congruence Relation}.
    $$
    [X]=[Y]
    $$
whenever $X$ and $Y$ are nonsingular projective varieties that
give the {\it same virtual Chow motive}.  We will uncoil this
definition a bit below.  All that is `motivic' about motivic
measure stems from this particular congruence relation.

\begin{definition}
The quotient of the free abelian group by the scissor and
congruence relations is the motivic scissor ring $\ring{K}$. (The
letter `K' is the standard notation for a Grothendieck group,
which for our purposes is just another name for a scissor group.)
The localized version $\ring{K}[\ring{L}^{-1}]\otimes \ring{Q}$
will be called the {\it localized motivic scissor ring} and
denoted $\S{mot}$. (It will become clear in
Section~\ref{sec:scale} why it is useful to invert $\ring{L}$.)
\end{definition}

It is time to uncoil the definition of this congruence relation.
There is a category of Chow motives.  To describe this category,
we assume familiarity with the Chow groups $A^i(X)$ of a variety
$X$.  They are groups of cycles of a given codimension $i$ modulo
the subgroup of cycles that are rationally equivalent to $0$. A
detailed treatment of cycles, rational equivalence, and Chow
groups can be found in \cite{Ful}. Other good treatments of Chow
motives can be found in \cite{Sch} and \cite{VDG}.

An object in the category of Chow motives is a triple $(X,p,m)$
where $X$ is a smooth projective variety of dimension $d$, $p$ is
an element in the Chow ring $A^d(X\times X)$ that is a projector
($p^2=p$), and $m$ in an integer.  The set of morphisms from
$(X,p,m)$ to $(X,p',m')$ is defined to be the set
    $$
    p' A^{d+n-m}(X\times Y)p.
    $$
Varieties that are not isomorphic as varieties can very well
become isomorphic when viewed as Chow motives. For example,
isogenous elliptic curves are isomorphic as Chow motives.

There is a canonical morphism from the Grothendieck ring of the
category $\op{Var}_{\ring{Q}}$ to the Grothendieck ring of the
category of Chow motives.  We let $\ring{K}$ be the image of this
morphism.  To say that two varieties are equal as virtual Chow
motives is to say that they have the same class in $\ring{K}$.

\subsection{The motivic counting measure}

The following theorem follows from a deep investigation of Chow
motives, and the theory of quantifier elimination for
pseudo-finite fields.

\begin{theorem} \label{thm:ring-hom}
There exists a unique ring homomorphism $\S{cover}\to\S{mot}$ that
satisfies the following property (Zero Sets).
\end{theorem}

\noindent{\bf Zero Sets}.   If $\phi$ is a ring formula that is
given by the conjunction of polynomial equations, then $[\phi]$ is
sent to the affine variety defined by those polynomial equations.

\bigskip

There are ring homomorphisms $\S{count}\to\S{cover}\to\S{mot}$. We
use the notation $\phi\mapsto[\phi]$ for the class of $\phi$ in
any of these rings, depending on the context.

\begin{definition}
The composite map $\phi\mapsto[\phi]\in\S{mot}$ will be called the
{\it motivic counting measure\/} of the formula $\phi$.
\end{definition}

The motivic counting measure of a ring formula is thus represented
by a rational linear combination of varieties over $\ring{Q}$. I
like to think of the motivic counting measure as counting the
number of solutions of the ring formula over finite fields in a
way that does not depend on the finite field.  Instead of giving
the answer as a particular number, it gives the answer in terms of
a formal combination of varieties having the same number of
solutions over a finite field.  Here is the precise statement.

\begin{theorem}  Let $\phi$ be a ring formula, and let $\sum a_i
[X_i]$ be a representative of the motivic counting measure
$[\phi]$ as a formal linear combination of varieties.  Choose a
model of each $X_i$ over $\ring{Z}$.  For all $r$ and for all but
finitely many primes $p$, the number of solutions of $\phi$ in
$\ring{F}_{p^r}$ is equal to
    $$
    \sum a_i \#X_i(\ring{F}_{p^r}).
    $$
\end{theorem}

%\begin{example}  Let us work a trivial example.
%    Fix a nonzero integer $n$. Let $\phi_n$ be the ring formula that asserts
%    there exists a right triangle with area $n$.
%        $$
%        \phi_n : \exists a b c.\quad (2 n = a b) \land (a^2 + b^2
%        = c^2)\land (b\ne 0).
%        $$
%    The affine variety defined by the ideal $(2 n - a b,c^2-a^2-b^2,b b'-1)$ in
%    $\ring{A}^4$ is absolutely irreducible.  By the definition of
%    pseudo-finite field, this variety has a rational point in
%    every pseudo-finite field.  Thus, the sentence $\phi_n$ is
%    true in every pseudo-finite field.  The motivic counting
%    measure of $\phi_n$ is the class associated with a point
%    $[\op{pt}]$.
%\end{example}

\begin{example} As an example, let us calculate the motivic
counting measure of the `set' of nonzero cubes.  The formula is
given by
    $$
    \phi(x):\quad \lq\exists y.~(y^3 = x)~\land~(x\ne 0)\rq.
    $$
The scissor relation can be used to break $\phi$ into two disjoint
pieces $\phi =\phi_1\lor\phi_2$: the part $\phi_1$ on which $-3$
is a square and the part $\phi_2$ on which it is not.   Let
$\ring{M}$ be the class in $\S{mot}$ corresponding to the
zero-dimensional variety $x^2+3=0$.  The class $\ring{M}$ has two
solutions or no solutions according as $-3$ is a square or not.
When $-3$ is a square, the cube roots of unity lie in the field,
so that the nonzero points on the affine line give a $3$-fold
cover of $\phi_1$ (under $y\mapsto y^3$).  Thus, $\phi_1$ has
measure
    $$
    \left(\frac{\ring{L}-1}{3}\right)\frac{\ring{M}}{2}.
    $$
On the other hand, if $-3$ is not a square, each non-zero element
of a pseudo-finite field of characteristic zero is a cube, so that
$\phi_2$ has measure
    $$
    ({\ring{L}-1})\left(1-\frac{\ring{M}}{2}\right).
    $$
The sum of these two terms is the measure of $\phi$ in $\S{mot}$.
\end{example}

%\subsubsection{What are Galois formulas?}%

%If $X$ is closed in $\ring{A}^n$, then a morphism $f:Y\to X$
%determines a ring formula in the free variables
%$x=(x_1,\ldots,x_n)$.
%    $$
%    \phi_f(x) :\quad \exists y \in Y.\quad (f(y) = x)\land (x\in
%    X).
%    $$
%(With a little work, this condition can be expanded into a ring
%formula.)

%Now assume that $f:Y\to X$ is a Galois cover with Galois group
%$G$.  Let $C$ be an abelian subgroup of $G$.  Then there is a ring
%formula $\phi_{Y,X,C}(x)$ that asserts that $x\in\ring{A}^n$ is
%the image of $y\in Y/C$ and is not the image of any $y\in Y/C'$
%for any proper subgroup $C'\subset C$.

%Equation~\ref{eqn:galois} has a nice geometric interpretation. Let
%$K$ be any pseudo-finite field.   The map from $Y/C$ to $X$ sends
%the set of rational solutions of $\phi_{Y,Y/C,C}$ onto the set of
%rational solutions of $\phi_{Y,X,C}$.  This is a $N_G(C)/C$ to $1$
%map.  This relation may be viewed as a special case of a
%similarity relation (as opposed to a congruence relation) in the
%scissor ring: if one formula $\phi$ is an $n$-fold ``cover'' of
%another formula $\phi'$, then $[\phi] = n [\phi']$.

\section{Locally Compact Fields and Haar Measures}

This section makes the transition from finite fields to locally
compact fields and from counting measures to additive Haar
measures.

In Section~\ref{sec:counting}, we developed a universal counting
measure for ring formula.  It may be viewed as counting solutions
to the ring formula over a finite field in a way that does not
depend on the finite field.

Counting measures are a rather simple and uninteresting type of
measure.  In this section, we construct a universal (motivic)
measure with ties to locally compact fields. This new measure may
be viewed as the volume expressed in a way that does not depend on
the locally compact field.  To carry out the construction, we must
work with a different collection of formulas (called DVR formulas)
that are better adapted to locally compact fields.  `DVR' is an
acronym for {\it discrete valuation ring}.

\subsection{Examples of rings}

To make the transition from finite fields to locally compact
fields, we wish to replace ring formulas with formulas in a
language that has a rich assortment of locally compact structures.

\begin{example}  Let $\ring{C}[[t]]$ be the ring of formal power
series with complex coefficients.  A typical element of this ring
has the form
    $$
    x = \sum_{i=k}^\infty a_i t^i
    $$
(with no constraints on the convergence of the series).  Pick the
initial index $k$ so that $a_k\ne 0$ (if $x\ne 0$).

The {\it valuation\/} of $x$ is defined to be the integer $k$:
    $$
    \op{val}(x) = k.
    $$
The {\it angular component\/} of $x$ is defined to be the complex
number $a_k$.
    $$
    \op{ac}(x) = a_k\in\ring{C}^\times.
    $$
(In the special case $x=0$, we set $\op{val}(0)=\infty$ and
$\op{ac}(0) = 0$.)
\end{example}

The name {\it angular component\/} is not meant to suggest any
precise connection to angles.   The name is based on a loose
analogy with the polar coordinate representation of a complex
number: just as the angular component $\theta$ of a nonzero
complex number $re^{i\theta}$ distinguishes among complex numbers
of the same magnitude (or valuation) $r$,  so the angular
component of a formal power series helps to distinguish among
formal power series of a given valuation $k$.

There are many other rings with similar functions, $\op{ac}$ and
$\op{val}$.  For example, we can change the coefficient ring of
the formal power series from $\ring{C}$ to any other field $k$ to
obtain $k[[t]]$.  Or we can take the field of fractions of
$k[[t]]$, which is the field of formal Laurent series with
coefficients in $k$:
    $$
    k((t)) = \{\sum_{-N}^\infty a_i t^i ~|~ a_i \in K\}.
    $$

For each prime $p$, there are valuation and angular component
functions defined on the field of rational numbers.  If $x$ is a
nonzero rational number, pick integers $a,b,c,N$ so that
    $$
    x = a p^N + \frac{b p^{N+1}}{c},
    $$
where $c$ is not divisible by $p$, and $a \in\{1,\ldots,p-1\}$.
The integers $a$ and $N$ are uniquely determined by this
condition. Define the valuation of $x$ to be
$\op{val}_p(x)=N\in\ring{Z}$ and the angular component of $x$ to
be image of $a$ modulo $p$ in $\ring{F}_p$.

\begin{example}
    If $p=2$ and $x=17/8$, then
        $$
        17/8 = 1.2^{-3} +2,\quad \op{val}_2(17/8) = -3,\quad
        \op{ac}(17/8) = 1\in\ring{F}_2.
        $$
\end{example}

Other examples, can be obtained from this one by completion.  For
each $p$,
    $$
    d(x,y) = (1/2)^{\op{val}_p(x-y)}
    $$
is a metric on the set of rational numbers.    The completion is a
locally compact field, called the field of $p$-adic numbers
$\ring{Q}_p$. The valuation $\op{val}_p$ and angular component
function $\op{ac}$ functions extend to the completion.

\subsection{The DVR language}

We have seen by example that there are many rings with functions
$\op{val}$ and $\op{ac}$.   In each case, there are three separate
rings that come into play: the domain of the functions $\op{val}$
and $\op{ac}$, the range of the function $\op{val}$ (which we
augment with a special symbol $\{\infty\}$ for the valuation of
$0$), and the range of the function $\op{ac}$. We call these rings
the valued ring, the value group, and the residue field,
respectively.

We formalize this relationship as a language in first-order logic
with function symbols $\op{val}$ and $\op{ac}$.  We allow
ourselves to build syntactically well-formed expressions with
variables, parentheses, quantifiers, the function symbols
$\op{val}$ and $\op{ac}$, the usual ring operations
$(0,1,(+),(-),(*),(=))$ on the valued ring and residue field, and
the usual group operations and inequalities on the value group
$(0,(+),(\le))$. These formulas will use variables of three
different types $x_i$ for the value ring, $m_i$ for the value
group, and $\xi_i$ for the residue field. Quantifiers $\forall$,
$\exists$ can be used to bind all three sorts of variables.

The construction of first-order languages is commonplace in logic,
but even without any background in logic, it is not hard to guess
whether a formula is syntactically correct.  We allow standard
mathematical abbreviations similar to those introduced above for
ring formulas.
    $$
    \lq\forall y.~(\exists x.~ x^2=y) \implies (\exists m.~2m =
    \op{val}(y)).\rq
    $$
is syntactically correct.  But
    $$
    \lq\forall f.~\forall x.~\forall y. f(y,\op{ac}(y))\rq
    $$
is not well-formed, because quantifiers are not allowed over
higher-order relations $f$ in a first-order language.  Also,
    $$
    \lq\forall x~\xi.~ (0\le x) \lor (\op{ac}(x) = \xi)\rq
    $$
is not well-formed, because of a type error; the variable symbol
$x$ appears once as an integer $0\le x$ and again as variable in
the valued field $\op{ac}(x)$.

A syntactically correct formula is called a DVR formula. {\it The
aim of motivic measure is to compute the ``volume'' of a DVR
formula in a universal way; that is, in a way that does not depend
on the underlying locally compact field.}

\subsection{Assumptions on the ring}

The various examples that we have mentioned are all structures for
the DVR language: rings of formal power series $k[[t]]$, fields of
formal Laurent series $k((t))$. For each prime $p$,
$(\ring{Q},\op{ac},\op{val}_p)$ is a structure for the language,
as well as its completion $(\ring{Q}_p,\op{ac},\op{val}_p)$.

We will temporarily restrict the set of examples to structures
$(K,k,\op{ac},\op{val})$ that satisfy the following conditions.
    \begin{itemize}
        \item $K$ is a valued field of characteristic zero, with
        valuation function $\op{val}:K\to\ring{Z}\cup\{\infty\}$
        and angular component
        functions $\op{ac}:K\to k$.
        \item The residue field $k$ has characteristic zero.
        \item $K$ is henselian.  (We review the definition below.)
        \label{item:hensel}
    \end{itemize}
Examples that satisfy these conditions include the fields
$k((t))$, where $k$ has characteristic zero.  The analogy that
will guides us is that these fields stand in the same relation to
locally compact DVR fields, as pseudo-finite fields do to finite
fields.

\subsection{Henselian field}

There is only one plausible definition for a henselian field:  A
field is {\it henselian} if the field satisfies Hensel's lemma.

Hensel's lemma gives checkable conditions on a polynomial that
insure that it has a root in a given neighborhood.   Hensel's
lemma occupies same ground in the realm of DVR rings that the
intermediate value theorem occupies in the realm of real numbers.
(The intermediate value theorem also gives checkable conditions on
a polynomial that insure that it has a real root in a given
neighborhood.)

Our experience with motivic counting measures has alerted us to
the importance of quantifier elimination, that is, the process of
replacing a formula with quantifiers $\forall,\exists$ with an
equivalent formula that does not contain quantifiers.  The
simplest case of quantifier elimination is the determination of
when there exists a root of a polynomial. Without a criterion for
the existence of roots to polynomials, quantifier elimination
would be impossible. For the pseudo-finite fields, this is handled
through the defining property of pseudo-finite fields that ``every
absolutely irreducible variety has a root.'' For real fields,
quantifier elimination is based on the intermediate value theorem.
For henselian fields, quantifier elimination is based on Hensel's
lemma.

\begin{lemma} (Hensel's lemma) For every monic polynomial $f\in
K[x]$, whose coefficients have non-negative valuation, and for
every $x$ such that
    $$\op{val}(f(x)) >0$$
and
    $$\op{val}(f'(x)) =0,$$
there exists $y\in K$ such that $f(y)=0$ and $\op{val}(y-x)>0$.
\end{lemma}

This is stated as a lemma, but we view it as a condition on the
field $K$ and its valuation.  It can be proved that the fields
$k((t))$ and $\ring{Q}_p$ are henselian by showing that under the
hypotheses of Hensel's lemma, Newton's approximations to the roots
    $$
    \begin{array}{lll}
    x_0 &= x\\
    x_{n+1} &= x_n - f(x_n)/f'(x_n)
    \end{array}
    $$
converge to a root.

\subsection{Quantifier elimination}

\begin{theorem} (Pas \cite{pas})  Let $K$ be a field satisfying
the other conditions enumerated in \ref{item:hensel} with residue
field $k$ . Let $\phi$ be a DVR formula.   Then there is another
formula $\phi'$ without quantifiers of the valued field sort such
that
    $$
    \forall (x,\xi,m)\in K^n \times k^m\times
    (\ring{Z}\cup\{\infty\})^r.\quad \phi^K(x,\xi,m) =
    \phi^{\prime~K}(x,\xi,m).
    $$
Moreover, the formula $\phi'$ can be chosen to be independent of
the structure $K$.
\end{theorem}

\subsection{Outer measure of a DVR formula}

As a first step toward constructing the measure of a DVR formula,
we will define an outer measure of a formula.  To motivate this
construction, it might be helpful first to describe an analogous
construction in Euclidean space.

\subsubsection{An outer measure in Euclidean space}

Fix a positive integer $m$.  Tile Euclidean space with cubes of
width $1/2^m$ whose vertices are centered at points $a$ with
coordinates $a_i\in \ring{Z}/2^m$.

According to the Calculus~$101$ approach to volume, we can
approximate the volume of a set by counting the number of cubes
that it meets. Let $A$ be a bounded set in $\ring{R}^n$. Let
$C_m(A)$ be the set of cubes in this tiling that meet $A$.  In our
naive approach to measure, let us define the outer measure of $A$
at level $m$ in dimension $n$ to be
    \begin{equation}
    \frac{\#C_m(A)}{2^{mn}},
    \label{eqn:cubes}
    \end{equation}
that is the number of cubes divided by the scaling factor
$2^{mn}$.   (If doing so did not involve logical circularity, we
would identify $1/2^{mn}$ with the volume of cube and the entire
expression as the volume of the set $C_m(A)$ of cubes.)

\begin{figure}[htb]
  \centering
  \includegraphics{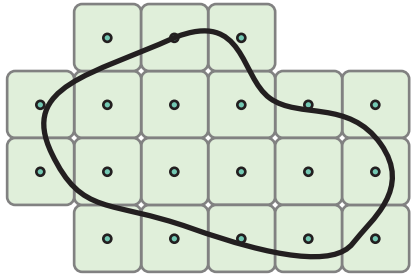}
  \caption{Volumes of DVR formulas can be approximated in Calculus
  101 fashion by counting centers of cubes
   that meet a given formula, scaled according to the size of the cubes.}
  \label{fig:cubes}
\end{figure}

The outer motivic measure of a DVR formula will be formed in an
entirely analogous way.  Of course, we will need to decide what to
use for cubes, how to count the number of cubes that ``meet'' a
given formula, and what scaling factor to use.  Once we make these
decisions, the formula for outer measure will take precisely the
same form as Equation~\ref{eqn:cubes}.

In the planar case, we gave a construction of area of polygons as
taking values in a scissor group $\S{poly}$. The outer
approximation of any bounded planar set $A$ by squares gives a
value in the scissor group of polygons.  Here too, if our outer
approximation to a DVR formula is with a ring formula, then the
value of the outer measure of the DVR formula will be in a scissor
ring $\S{mot}$.

Given all our preliminaries, it almost goes without saying that
the number of cubes appearing in the numerator of
Equation~\ref{eqn:cubes} will be replaced with the motivic
counting measure of a ring formula.

\subsubsection{Cubes}

What is a cube?  Well, it is a product of equal width intervals.
In DVR formulas, a cube centered at $a$ of ``width'' $m$ is again
a product of intervals:
    $$
    \{(x_1,\ldots,x_n)\in K^n ~|~ \op{val}(x_i-a_i) \ge m,\quad\text{ for }i=1,\ldots,n\}.
    $$
If $K=k[[t]]$, then the interval around a formal power series $a$
is the set of all formal power series with the same leading terms.
Shaking (wagging) the tails of the power series fills out the
interval. In other words, we can make precise the idea of covering
a DVR formula with cubes by replacing each solution to the DVR
formula with a bigger set where the tails of the solutions are
allowed to vary.

Let us make this precise. We have truncation map
    $$
    \begin{array}{lll}
        k[[t]]\to k[[t]]/(t^m) \simeq k^m\\
        \\
        \sum_{0}^\infty a_i t^i \mapsto \sum_{0}^{m-1} a_i t^i
            \mapsto (a_0,\ldots,a_{m-1}).
    \end{array}
    $$
In the opposite direction, given $b\in K^m$, there is  a
polynomial with those coefficients
    $$
    p(b,t) = \sum_{0}^{m-1} b_i t^i \in k[[t]]
    $$

\begin{definition} Let $\phi$ be a DVR formula with free variables $(x_1,\ldots,x_n)$
and no free variables of other sorts.   An {\it outer ring
formula} $\phi_m$ approximation to $\phi$ {\it (at level $m$)} is
a ring formula in $n m$ free variables $u_{ij}$ such that over
every field $k$:
    $$
    \begin{array}{lll}
        &\{u \in k^{n m}~|~ \phi_m(u)\} =\\
    &\quad\quad\{u \in k^{n m}~|~\exists a_1,\ldots,a_n.~
    \phi(a_1,\ldots,a_n)~\land~ \op{val}(a_i-p(u_{ij},t))\ge m\}.
    \end{array}
    $$
\end{definition}
This set is the set of centers of cubes that contain a solution to
$\phi$.

\begin{theorem}
Outer ring formula approximations exist for every DVR formula
$\phi$ at every level $m$.
\end{theorem}

The proof of this theorem uses quantifier elimination results to
eliminate the quantifiers that bind variables ranging over the
valued field. It uses results of Presburger on quantifier
elimination to eliminate the quantifiers that range over the
additive group of integers.  The quantifiers that bind variables
in the residue field remain as quantifiers in the ring formula
$\phi_m$.

\subsubsection{Scaling Factors}\label{sec:scale}

How is the scaling factor chosen in Equation~\ref{eqn:cubes} for
Euclidean outer measures?  The scaling factor $1/2^{n m}$ is the
unique constant that has the property that if the set $A$ is
itself a union of properly aligned cubes (of width $m'$), then the
outer measure of $A$ is independent of $m$ for all $m\ge m'$.

To find the scaling factor for DVR formulas, we work a simple
example in which the DVR formula is itself a union of cubes of
width $m'$ (that is, its set of solutions is stable under
perturbation of the power series tails).

\begin{example}  Let $\phi(x_1,\ldots,x_n) = \ring{T}$, a formula
that is true for all values of the free variables $x_i$.  In this
case the outer ring formula approximation is exact.  Substitute
polynomials $p(u_{i\cdot},t)$ for each $x_i$ and expand in terms
of $m n$ distinct free variables $u_{ij}$ to get
    $$
    \phi_m(u_{ij}) = \ring{T}
    $$
for all input values $u_{ij}$.  The number of solutions of
$\phi_m$ over a finite field $\ring{F}_q$ is $q^{nm}$. If we take
the motivic counting measure of $\phi_m$, we find that the variety
that counts the points of $\phi_m$ over any finite field is the
affine space of dimension $n m$:
    $$
    \# \ring{A}^{n m}(\ring{F}_q) = q^{n m}.
    $$
The class of $\phi_m$ in $\ring{K}[\ring{L}^{-1}]\otimes \ring{Q}$
is
    $$[\ring{A}^{n m}] = [\ring{A}^1]^{n m} = \ring{L}^{n m}.
    $$
\end{example}

From this one example, we see that the scaling factor for DVR
formulas must be $1/\ring{L}^{n m}$.

\begin{definition} Let $\phi$ be a DVR formula.  Let the {\it outer
measure\/} of $\phi$ at level $m$ be given by
    $$
    \frac{[\phi_m]}{\ring{L}^{n m}}\in
    \ring{K}[\ring{L}^{-1}]\otimes\ring{Q}=\S{mot}.
    $$
\end{definition}

This formula is analogous to Formula~\ref{eqn:cubes} for the
Euclidean outer measure at level $m$.  The numerator counts the
number of centers of cubes that contain a solution to the DVR
formula.

\begin{definition} Let the {\it motivic measure (or motivic volume)} of
$\phi$ be given by
    $$
    \lim_{m\to\infty} [\phi_m]\ring{L}^{-n m},
    $$
whenever that limit exists. (The limit must be taken in a
completion of $\S{mot}$.)
\end{definition}

\subsection{The universal nature of motivic measure}

Just as the motivic counting measure counts solutions to ring
formulas over finite fields in a field independent way, so the
motivic measure takes the volume of a DVR formula over locally
compact fields in a field independent way.\footnote{It is
impossible for the structure $K$ both to be locally compact and to
have a residue field $k$ of characteristic zero, as required by
Condition~\ref{item:hensel}. The residue field of a locally
compact field is always finite.  In these final paragraphs, we
allow the residue field to have positive characteristic.  }

There is a good theory of measure  on locally compact fields. This
is the Haar measure, which is  translation invariant. Given a DVR
formula $\phi$ and a locally compact structure $K$ with ring of
integers $O_K$, we can take the volume of the set of solutions to
the DVR formula
    \begin{equation}
    \op{vol}(\{x\in O_K^n ~|~ \phi^K(x)\},dx).
    \end{equation}
The measure $dx$ can be given a canonical normalization by
requiring that it assigns volume $1$ to the full set $O_K^n$.

We are now ready to state the main result on motivic measure. Like
all the other principal results in this article, the result is due
to J. Denef and F. Loeser.

\begin{theorem}\label{thm:main}
The motivic volume of $\phi$ is universal in the following sense.
Let $\sum a_i [X_i]\ring{L}^{-N_i}$ be any representative of the
motivic volume of $\phi$ as a convergent formal sum of varieties
over $\ring{Q}$. Pick models for the varieties over $\ring{Z}$.
After discarding finitely many primes, for any locally compact
structure of the DVR language, the $K$-volume of the formula is
given by a convergent sum (in $\ring{R}$)
    $$
    \sum a_i \#X(\ring{F}_q) q^{-N_i},
    $$
where $\ring{F}_q$ is the residue field of $K$.
\end{theorem}

This wonderful result states that the Haar measures on all locally
compact fields have an deep underlying unity.  The volumes of sets
can be expressed geometrically in a way that is independent of the
underlying field.

Moreover, there are effective procedures to calculate the
varieties $X_i$ and the coefficients $a_i,N_i$ that represent the
outer motivic volume at level $m$. If the outer ring formula
approximations $\phi_m$ converge at some finite level $m$ to the
DVR formula $\phi$, then we obtain effective procedures to
calculate the motivic volume of the formula.

\section{Applications and Conclusions}

What good is motivic measure?  Here are a few examples.

\subsection{Invariants of ring formulas} The group $\S{mot}$ is
generated by varieties $\op{Var}_{\ring{Q}}$.  Many geometrical
invariants of varieties (such as Euler characteristics and Hodge
polynomials) can be reformulated as invariants of the ring
$\S{mot}$.   This gives a novel way to attach invariants to every
ring formula $\phi$: take a geometric invariant of
$[\phi]\in\S{mot}$.  In particular, ring formulas have Euler
characteristics and Hodge polynomials! For example, the formula
for the nonzero squares in a field
    $$
    \lq\exists y.~ (y^2 = x) ~\land~(x \ne 0)\rq
    $$
has Euler characteristic zero.

\subsection{Geometry of varieties.}

There is a motivic change-of-variables formula that is similar to
the standard change of variables formula in calculus.  Using this
formula, it is sometimes possible to show that two birationally
equivalent varieties have the same motivic volume.  This has deep
implications for the geometry of the two varieties.  In
particular, the motivic volume determines the Hodge polynomial of
the varieties.

This approach was followed by Kontsevich, who used a
change-of-variables calculation to show that birationally
equivalent projective Calabi-Yau manifolds have the same Hodge
numbers \cite{kont}.  Applications to orbifolds appear in
\cite{luper}.

\subsection{Computation of $p$-adic integrals}

Many integrals over $p$-adic fields are notoriously difficult to
calculate.  Motivic measure exposes the underlying similarities
between volumes on different $p$-adic fields.  It gives a decision
procedure to calculate $p$-adic integrals (at least when the data
defining the integral can be expressed as DVR formulas that can be
reproduced at some finite level $m$).   In particular, this means
that a computer can be programmed to compute a large class of
$p$-adic integrals.

\subsection{Generating Functions}

Motivic counting gives a way of counting that is independent of
the finite field.   Let
    $$
    Z_p(t) = \sum_{i=0}^\infty a^{(p)}_i t^i
    $$
be a generating function, where the constants $a^{(p)}_i$ are
obtained by counting solutions to a formula in some $p$-dependent
way. (Each generating function depends on a single prime $p$.)
Motivic measure can often give a way of forming a $p$-independent
series
    $$
    Z_{\op{mot}}(t) = \sum_{i=0}^\infty [a_i] t^i
    $$
taking values in $\S{*}[[t]]$ and specializing for almost all $p$
to the $p$-dependent series $Z_p(t)$.  The motivic series collects
the behavior of the various series $Z_p(t)$ into a single series.

Denef and Loeser have studied motivic versions of Hasse-Weil
series, Igusa series, and Serre series.  They have used the
general motivic series to prove that various properties of these
series are independent of the prime $p$.  See \cite{DLR}.

\subsection{Concluding Remarks}

This article is an exposition of a particular version of motivic
integration, called {\it arithmetic motivic integration}. Proofs
of results stated in this article can be found in \cite{ICM} and
\cite{DL}.  Motivic integration has been developing at a
break-neck pace, ever since Kontsevich gave the first lecture on
the topic in 1995.  The version of motivic integration developed
in the late nineties goes by the name of {\it geometric motivic
integration}. Geometric motivic integration is a coarser theory,
but is sufficient for many applications.  Good introductions are
\cite{craw} and \cite{L}.  Some articles on geometric motivic
integration include \cite{DLG} and \cite{DLM}. Another version of
motivic integration has been developed by J. Sebag for formal
schemes \cite{sebag}. See also \cite{LS}. Cluckers and Loeser are
in the final stages of preparation of an ultimate version of
motivic integration that subsumes both geometric and arithmetic
motivic integration \cite{CL}.

We began this article by stating that motivic measure does not fit
neatly into the tradition of Hausdorff, Haar, and Lebesgue.
However, a major result states that the motivic measure
specializes to the additive Haar measure on locally compact fields
(Theorem~\ref{thm:main}).  Thus, the motivic measure is perhaps
not so peculiar after all.  In fact, in many respects it is
strikingly similar to the additive Haar measure on locally compact
fields.  It has been my experience when I calculate motivic
volumes to lose track -- mid-calculation -- of which measure is
being used.

\section*{Appendix: Historical Remarks on Galois Stratification}

 \medskip
 \centerline{by Michael Fried}
 \bigskip

Fran\c{c}ois Loeser was one of the (three, including myself and
Gross) foreigners who gave talks at Yasutaka Ihara's sixtieth
Birthday Conference several years ago in Tokyo. I was totally
unaware of the relation between Denef and Loeser's work and mine,
until he gave that talk. My story starts with my 1976 Annals paper
\cite{FAn}.

This featured Galois stratification: a dissection of Diophantine
statements to produce uniform computations for rational points
running over primes. It also worked over all extensions of a finite
field and through $p$-adic integration over all $p$-adic
completions. My followup paper \cite{FL} (written in 1978 from
lectures I gave at Yale), measured the bad primes on general
statements. I detected bad primes through failure of the Euler
factor for that prime to fit the uniform variation of Euler factors.

Galois stratification takes up considerable space in the papers of
Denef and Loeser, for the simple reason that quantifier elimination
is the tool that allows the conclusions you taut in
Section~\ref{sec:qe}.

General diophantine statements to me included such problems as came
from Artin's conjecture: Failure of a hypersurface of degree $d$ in
projective $d^2$ space to have a $p$-adic point. You can see the
territory was that defined by Ax and Kochen. (I went with Ax to
Stony Brook, instead of going with a tenure offer to University of
Chicago.) You also see it is a more precise problem than considered
by the motivic counting of points on varieties over finite fields.
That considers only varieties with coefficients in $\ring{Q}$
(blindly throwing out any knowledge a finite set of bad primes).
There weren't any motives when I did my paper, only the category of
Galois stratifications that I invented.

Loeser's first words -- I've seen him since, but never had talked to
him before -- to me, at the train station back to our hotel in Tokyo
were: ``How come nobody knows about this?" He was holding up the
Fried-Jarden book opened to Chap. 25. That was friendly and generous
to me. Also, deserved! It is the only gesture in my direction made
by anyone in the last 25 years toward me on the topic.

I've made little noise about the neglect of my early papers. It came
from an antagonism
%by people like Lang -- not to me personally, I
%think because I'd never met him, but to
toward the ``school" from which I came. Nevertheless, it killed my
career effectively for 12 years, until John Thompson asked me to
renovate my Inverse Galois Problem ideas.

The technical heart of Denef and Loeser's results are Galois
stratification, in the original form conceived in my Annals paper.
Their goals were versions of mine in \cite{FL}.

Did I have the idea of using an abstract gadget (Galois
stratification) to ``measure'' the truth of statements over
arithmetic rings. Yes, I did in the early seventies! Did I
acknowledge the motivic approach, as soon as I heard of it (many
years later; especially the Denef-Loeser idea of Hodge invariants)?
Yes! Did my approach, technique and persistence influence the topic?
Yes!
%
%Here's a different type of question.
Do I now have any influence over the subject%, or can one safely
%dismiss my part in it?
? Not clear!
%Yet, I just refereed (assiduously and with perspicacity) then
%accepted Nicaise's paper ``Relative motives and the theory of
%Pseudo-finite fields" for IRMN.

%True, my work in the last 20 years doesn't concentrate on topics
%related to these ideas. So, one might question whether there are
%still any contributions left from my approach. The enclosed paper
%"The place of exceptional covers in all diophantine relations" is a
%project to notice that certain special formulas have a significant
%place among all possible formulas. Exceptional covers and their
%generalization, Davenport pairs, give those formulas. These produce
%what I term as ``universal relations" between Poincar\'e series.
%This will be the main theme of the followup paper.

    \renewcommand{\thefootnote}{}
    \footnote{Version: May 5, 2005.}
    \footnote{Work supported by the NSF}
    \footnote{Copyright (c) 2004, 2005, Thomas C. Hales.}
    \footnote{This work is licensed under the Creative Commons Attribution
License. To view a copy of this license, visit
http://creativecommons.org/licenses/by/1.0/ or send a letter to
Creative Commons, 559 Nathan Abbott Way, Stanford, California
94305, USA.}

\end{document}